\documentclass[12pt,dvipdfmx]{article}
\usepackage{amsmath}
\usepackage{amssymb,amscd}
\usepackage{bm}
\usepackage{wrapfig}
\usepackage{color}
\usepackage{bm}
\usepackage[all]{xy}

\usepackage{theorem}

\newtheorem{theorem}{Theorem}[section]
 \newtheorem{proposition}[theorem]{Proposition}
 
 \newtheorem{lemma}[theorem]{Lemma}

{\theorembodyfont{\normalfont}
 \newtheorem{definition}[theorem]{Definition}
 \newtheorem{example}{Example}[section]
 \newtheorem{remark}{Remark}[section]
 }

\newcommand{\qed}{\nobreak \ifvmode \relax \else
      \ifdim\lastskip<1.5em \hskip-\lastskip
      \hskip1.5em plus0em minus0.5em \fi \nobreak
      \vrule height0.75em width0.5em depth0.25em\fi}


\newlength{\minitwocolumn}
\usepackage{amsmath}
\usepackage{amssymb,amscd}
\usepackage{bm}
\usepackage{amsbsy} 
\usepackage{wrapfig}
\usepackage{theorem}

\newcommand{\beq}{\begin{equation}}
\newcommand{\eeq}{\end{equation}}
\newcommand{\bea}{\begin{eqnarray*}}
\newcommand{\eea}{\end{eqnarray*}}
\newcommand{\beqa}{\begin{eqnarray}}
\newcommand{\eeqa}{\end{eqnarray}}

\def\bR{{\mathbb{R}}}

\newcommand{\calL}{{\mathcal L}}
\newcommand{\calM}{{\mathcal M}}

\newcommand{\calP}{{\mathcal P}}
\newcommand{\calS}{{\mathcal S}}

\newcommand{\bracket}[2]{\langle #1,\,#2\rangle}

\newcommand{\rd}{\mathrm{d}}

\newcommand{\uu}{{\underline{u}}{}}
\newcommand{\uv}{{\underline{v}}{}}
\newcommand{\ue}{{\underline{e}}{}}
\newcommand{\ualpha}{{\underline{\alpha}}{}}

\newcommand{\momega}{{\omega}}
\newcommand{\tmu}{{\tilde{\mu}}}
\newcommand{\hmu}{{\widehat{\mu}}}
\newcommand{\tH}{{\tilde{H}}}

\begin{document}


\baselineskip 0.7cm

\begin{titlepage}
\begin{flushright}
\end{flushright}

\vskip 1.35cm
\begin{center}
{\Large \bf
Homotopy momentum sections on multisymplectic manifolds
}
\vskip 1.2cm
Yuji Hirota${}^a$
\footnote{E-mail:\
hirotaATazabu-u.ac.jp
}
and
Noriaki Ikeda${}^b$
\footnote{E-mail:\
nikedaATse.ritsumei.ac.jp
}
\vskip 0.4cm
{\it
${}^a$School of Veterinary Medicine, 
Azabu University, \\
Sagamihara, Kanagawa 252-5201, Japan
}
\vskip 0.4cm

{\it
${}^b$Department of Mathematical Sciences,
Ritsumeikan University \\
Kusatsu, Shiga 525-8577, Japan \\
}
\vskip 0.4cm

\today

\vskip 1.5cm

\begin{abstract}
We introduce a notion of a homotopy momentum section 
on a Lie algebroid over a pre-multisymplectic manifold.
A homotopy momentum section is a generalization of the momentum map with a Lie group action and the momentum section on a pre-symplectic manifold, 
and is also a generalization of the homotopy momentum map on a multisymplectic manifold.
We show that a gauged nonlinear sigma model with Wess-Zumino term with Lie algebroid gauging has the homotopy momentum section structure.
\end{abstract}
\end{center}
\end{titlepage}

\tableofcontents

\setcounter{page}{2}


\rm

\section{Introduction}
\noindent
The moment(um) map is a fundamental object in symplectic geometry
inspired by conserved quantities and symmetries in the analytical mechanics.
In many cases, for a Lie group action on a symplectic manifold,
there exists a corresponding function called a Hamiltonian function.
Such a space is called a Hamiltonian $G$-space and the geometric
structure is described by the momentum map theory.

In this paper, a generalization of a momentum map theory on a multisymplectic manifold is discussed.
There are many versions of a momentum map theory on a multisymplectic 
manifold, such as a momentum map differential form \cite{Carinena:1992rb, Gotay:1997eg}, a multimomentum map
\cite{Madsen:2010qp, Madsen:2011ru},
a homotopy momentum map \footnote{
It is called a \textit{homotopy} moment map in the paper \cite{Fregier:2013dda}.}
\cite{Fregier:2013dda},
a weak (homotopy) momentum map \cite{Herman:2017, Herman:2018box}, 
etc.
Some recent developments about moment maps on multisymplectic manifolds are, for instance, \cite{Ryvkin:2015pra, Ryvkin:2020}.

Recently, Blohmann and Weinstein \cite{Blohmann:2018} have proposed a generalization of a momentum map and a Hamiltonian $G$-space with a Lie algebra 
(Lie group) action to a Lie algebroid, inspired by the analysis of the Hamiltonian formalism of general relativity \cite{Blohmann:2010jd}.
A Lie algebroid \cite{Mackenzie} is a widely analyzed generalization of a Lie algebra.
A Lie algebroid is an infinitesimal object of a Lie groupoid analogous to the fact that a Lie algebra is an infinitesimal object of a Lie group.
The generalization of a momentum map is called a momentum section and 
the corresponding space is called a Hamiltonian Lie algebroid.
One of authors gave new examples in the constrained Hamiltonian mechanics and sigma models, and generalized a momentum section theory 
\cite{Ikeda:2019pef}. Moreover a momentum section has been generalized 
to Courant algebroid actions \cite{Ikeda:2021fjk}.

In this paper, we propose a generalization of momentum sections
to a pre-multisymplectic manifold and analyze its properties.
We call it a homotopy momentum section 
since it is also a generalization of a homotopy momentum map.
A Lie algebroid called a homotopy Hamiltonian Lie algebroid is introduced 
as a generalization of a Hamiltonian $G$-space.
Moreover a Lie algebroid generalization of a weak momentum map is defined.

We introduce a homotopy momentum section and weak homotopy momentum section
in Sections \ref{sec:HMS} and \ref{sec:WHMS}.
In Section \ref{sec:HMM}, relations of a homotopy momentum section on an action Lie algebroid with a homotopy momentum map and
a weak homotopy momentum map are discussed.
If $n=1$, i.e., a base manifold is a symplectic manifold 
and a Lie algebroid is the action Lie algebroid,
a homotopy momentum section is a momentum map, thus is also equivalent to
a homotopy momentum map.
If $n \geq 2$, on a multisymplectic manifold,
a homotopy momentum section is not necessarily a homotopy momentum map.
However, we prove that both are equivalent if a natural condition is satisfied 
on the action Lie algebroid.
It is emphasized that the condition, \textit{equivariance}, of (weak) homotopy momentum sections are important for comparisons.
We also discuss the relation between a weak homotopy momentum section and a weak homotopy momentum map on an action Lie algebroid.

In Section \ref{sec:GNLSM}, we show that a nonlinear sigma model with WZ 
term in $n$ dimensions has a homotopy momentum section structure 
if we consider gauging with the Lie algebroid gauge symmetry. 
The nonlinear sigma model is an important physical model which has many physical applications \cite{Wess:1971yu}.
One of the applications is the T-duality of string theory \cite{Buscher:1987sk, Buscher:1987qj}.
If we add the WZ term, the target space has a pre-multisymplectic structure 
additional to a Riemannian structure.
Gauging of the nonlinear sigma model with a Lie algebroid was analyzed in
\cite{Bouwknegt:2017xfi, Bugden:2018pzv, Chatzistavrakidis:2015lga, 
Chatzistavrakidis:2016jfz, Chatzistavrakidis:2016jci, 
Chatzistavrakidis:2017tpk, Ikeda:2019pef, Wright:2019pru}.
If we take the Hull-Spence ansatz for gauging of the nonlinear sigma model,
and suppose that the gauge symmetry is a Lie algebroid,
the gauge invariant condition of the gauged action functional requires
existence of a homotopy momentum section on the target pre-multisymplectic 
manifold.

This paper is organized as follows.
In Section 2, a Lie algebroid and related notion are prepared and
notation is fixed.
In Section 3, a homotopy momentum section is defined and some examples are discussed.
In Section 4, a weak homotopy momentum section and equivariance are introduced.
In Section 5, relations between a (weak) homotopy momentum section, 
and a (weak) homotopy momentum map are discussed.
In Section 6, the gauged nonlinear sigma model with WZ term is considered and a homotopy momentum section structure in this model is showed.
Section 7 is devoted to discussion and outlook.
In the Appendix, some formulas are summarized.

\section{Preliminary: Lie algebroid}
\noindent
In this section, we explain the background geometry of this paper.
Formulas of a Lie algebroid, differentials and connections are discussed.

\subsection{Lie algebroid}

\begin{definition}
Let $E$ be a vector bundle over a smooth manifold $M$.
A Lie algebroid $(E, \rho, [-,-])$ is a vector bundle $E$ with
a bundle map $\rho: E \rightarrow TM$ called the anchor map, 
and a Lie bracket
$[-,-]: \Gamma(E) \times \Gamma(E) \rightarrow \Gamma(E)$
satisfying the Leibniz rule,
\begin{eqnarray}
[e_1, fe_2] &=& f [e_1, e_2] + \rho(e_1) f \cdot e_2,
\end{eqnarray}
{where $e_i \in \Gamma(E)$ and $f \in C^{\infty}(M)$.}
\end{definition}
A Lie algebroid is a generalization of a Lie algebra and the space of vector fields on a smooth manifold.
\begin{example}[Lie algebra]
Let a manifold $M$ be one point $M = \{pt \}$. 
Then a Lie algebroid is a Lie algebra $\mathfrak{g}$.
\end{example}
\begin{example}[Vector fields]
If a vector bundle $E$ is a tangent bundle $TM$ and $\rho = \mathrm{id}$, 
then a bracket $[-,-]$ is a normal Lie bracket 
on the space of vector fields $\mathfrak{X}(M)$
and $(TM, \mathrm{id}, [-,-])$ is a Lie algebroid.
It is called a tangent Lie algebroid.
\end{example}

\begin{example}[Action Lie algebroid]\label{actionLA}
Suppose a smooth action of a Lie group $G$ to a smooth manifold $M$, 
{$M \times G \rightarrow M$.}
The differential of the map induces an infinitesimal action of the Lie algebra $\mathfrak{g}$ of $G$ on a manifold $M$.
Since $\mathfrak{g}$ acts as a differential operator on $M$,
the differential of the map
determines a bundle map $\rho: M \times \mathfrak{g} \rightarrow TM$.
Consistency of a Lie bracket requires that $\rho$ is 
a Lie algebra morphism such that
\begin{eqnarray}
~[\rho(e_1), \rho(e_2)] &=& \rho([e_1, e_2]),
\label{almostLA}
\end{eqnarray}
and we obtain a Lie algebroid 
$(E= M \times \mathfrak{g}, \rho, [-,-])$.
Here the bracket in left hand side 
of \eqref{almostLA} is a Lie bracket of vector fields.
This Lie algebroid is called an \textit{action Lie algebroid}.
\end{example}

\begin{example}[Poisson structure]\label{Poisson}
An important nontrivial Lie algebroid is a Lie algebroid induced from a Poisson structure. 

A bivector field $\pi \in \Gamma(\wedge^2 TM)$ is called a Poisson structure if $[\pi, \pi]_S =0$, where $[-,-]_S$ is a Schouten bracket on the space of multivector fields, $\Gamma(\wedge^{\bullet} TM)$.
For a Poisson bivector field $\pi$, a bundle map is defined as
$\pi^{\sharp}: T^*M \rightarrow TM$ by $\bracket{\pi^{\sharp}(\alpha)}{\beta}
= \pi(\alpha, \beta)$ for all $\beta \in \Omega^1(M)$.

Let $(M, \pi)$ be a Poisson manifold. Then, a Lie algebroid structure is induced on $T^*M$.
$\pi^{\sharp}: T^*M \rightarrow TM$ is the anchor map, and
a Lie bracket on $\Omega^1(M)$ is defined by the so called Koszul bracket,
\begin{eqnarray}
[\alpha, \beta]_{\pi} = \calL_{\pi^{\sharp} (\alpha)}\beta - \calL_{\pi^{\sharp} (\beta)} \alpha - \rd(\pi(\alpha, \beta)),
\end{eqnarray}
where $\alpha, \beta \in \Omega^1(M)$.
$(T^*M, -\pi^{\sharp}, [-, -]_{\pi})$ is a Lie algebroid.
\end{example}

\begin{example}[Twisted Poisson structure]\label{tPoisson}
If a bivector field $\pi \in \Gamma(\wedge^2 TM)$ and 
$H \in \Omega^3(M)$ satisfy
\begin{eqnarray}
&& \frac{1}{2}[\pi, \pi]_S 
= \bracket{\otimes^{3} \pi}{H},
\label{tPoisson1}
\\
&& \rd H =0,
\end{eqnarray}
$(M, \pi, H)$ is called a twisted Poisson manifold \cite{Klimcik:2001vg, Park:2000au, Severa:2001qm, Ikeda:2019czt}.
Here $\bracket{\otimes^{3} \pi}{H}$ is defined by 
$\bracket{\otimes^{3} \pi}{H}(\alpha_1, \alpha_2, \alpha_3)
:= H(\pi^{\sharp} (\alpha_1), \pi^{\sharp} (\alpha_2), \pi^{\sharp} (\alpha_3))
$ for $\alpha_i \in \Omega^1(M)$.

If we take the same bundle map, $\pi^{\sharp}: T^*M \rightarrow TM$ 
as in Example \ref{Poisson},
and a Lie bracket on $\Omega^1(M)$,
{
\begin{eqnarray}
[\alpha, \beta]_{\pi,H} = \calL_{\pi^{\sharp} (\alpha)}\beta - \calL_{\pi^{\sharp} (\beta)} \alpha - \rd(\pi(\alpha, \beta))
+ \iota_{\pi^{\sharp}(\alpha)} \iota_{\pi^{\sharp}(\beta)} H,
\end{eqnarray}
}
for $\alpha, \beta \in \Omega^1(M)$.
Then, $(T^*M, -\pi^{\sharp}, [-, -]_{\pi, H})$ is a Lie algebroid.
\end{example}
One can refer to many other examples, for instance, in \cite{Mackenzie}.

\subsection{Lie algebroid differential and connection}\label{LAdiffandconn}
Consider the spaces of exterior products of sections of $E^*$ called 
the space of $E$-differential forms, $\Gamma(\wedge^{\bullet} E^*)$.
A \textit{Lie algebroid differential} ${}^E \rd: \Gamma(\wedge^m E^*)
\rightarrow \Gamma(\wedge^{m+1} E^*)$
such that $({}^E \rd)^2=0$ is defined as follows. 
\begin{definition}
A Lie algebroid differential ${}^E \rd: \Gamma(\wedge^m E^*)
\rightarrow \Gamma(\wedge^{m+1} E^*)$ is defined by
\begin{eqnarray}
{}^E \rd \alpha(e_1, \ldots, e_{m+1}) 
&=& \sum_{i=1}^{m+1} (-1)^{i-1} \rho(e_i) \alpha(e_1, \ldots, 
\check{e_i}, \ldots, e_{m+1})
\nonumber \\ && 
+ \sum_{1 \leq i < j \leq m+1} (-1)^{i+j} \alpha([e_i, e_j], e_1, \ldots, \check{e_i}, \ldots, \check{e_j}, \ldots, e_{m+1}),
\label{LAdifferential}
\end{eqnarray}
where $\alpha \in \Gamma(\wedge^m E^*)$ and $e_i \in \Gamma(E)$.
\footnote{In Equation \eqref{LAdifferential}, indices $i, j$ are not indices of local coordinates on $M$ but counting of elements of $\Gamma(E)$.}
\end{definition}

We introduce an ordinary connection on the vector bundle $E$,
$\nabla:\Gamma(E)\rightarrow \Gamma(T^*M \otimes E)$,
which is a 
$\bR$-linear map satisfying the Leibniz rule,
\begin{eqnarray}
\nabla (f e) = f \nabla e + (\rd f) \otimes e,
\end{eqnarray}
for $e \in \Gamma(E)$ and $f \in C^{\infty}(M)$.
A dual connection on $E^*$ is defined by the equation,
\begin{eqnarray}
\rd \bracket{\mu}{e} = \bracket{\nabla \mu}{e} + \bracket{\mu}{\nabla e},
\end{eqnarray}
for all sections $\mu \in \Gamma(E^*)$ and $e \in \Gamma(E)$.

The connection and the dual connection are extended to 
the \textit{exterior covariant derivative}
on the space of differential forms taking a value
in $E^{\otimes m} \otimes E^{* \otimes n}$, 
$\Omega^l(M, E^{\otimes m} \otimes E^{* \otimes n})
= \Gamma(M, \wedge^l T^*M \otimes E^{\otimes m} \otimes E^{* \otimes n})$.
The exterior covariant derivative is denoted by the same notation $\nabla$.
It increases the degree of the order of $\wedge^l T^*M$ by one.
i.e., $\nabla: \Omega^l(M, E^{\otimes m} \otimes E^{* \otimes n})
\rightarrow \Omega^{l+1}(M, E^{\otimes m} \otimes E^{* \otimes n})$.

On a Lie algebroid, another derivation is defined.
\begin{definition}
An \textit{$E$-connection} on a vector bundle $E'$ 
with respect to the Lie algebroid $E$ is a $\bR$-linear 
map
${}^E \nabla: \Gamma(E') \rightarrow {\Gamma(E^* \otimes E')}$ 
satisfying 
\begin{eqnarray}
{}^E \nabla_e (f e') = f {}^E \nabla_e e' + (\rho(e) f) e',
\end{eqnarray}
for $e \in \Gamma(E)$, $e' \in \Gamma(E')$ and $f \in C^{\infty}(M)$.
\end{definition}
The ordinary connection is regarded as an 
$E$-connection for $E=TM$, $\nabla = {}^{TM} \nabla$.

The \textit{standard $E$-connection} on $E$, 
${}^E \nabla: \Gamma(E) \rightarrow {\Gamma(E^* \otimes E)}$ is defined by
\begin{eqnarray}
&& {}^E \nabla_{e} e^{\prime} := \nabla_{\rho(e)} {e^{\prime}},
\label{stEconnection}
\end{eqnarray}
for $e, e^{\prime} \in \Gamma(E)$.
It is also denoted by 
$\rho^{\nabla}(e) e^{\prime} = {}^E \nabla_e {e^{\prime}}$
and called the covariantized anchor map.
For the tangent bundle $E'=TM$, there is the following canonical $E$-connection
on $TM$.
{If a normal connection $\nabla$ on $E$ as a vector bundle is given,}
a (canonical) $E$-connection on a tangent bundle, ${}^E \nabla: \Gamma(TM) \rightarrow {\Gamma(E^* \otimes TM)}$ is defined by
\begin{eqnarray}
{}^E \nabla_{e} v &:=& \calL_{\rho(e)} v + \rho(\nabla_v e)
= [\rho(e), v] + \rho(\nabla_v e),
\label{stEconnection1}
\end{eqnarray}
where 
$e \in \Gamma(E)$ and $v \in \mathfrak{X}(M)$.
It is also called the opposite connection.
For a $1$-form $\alpha \in \Omega^1(M)$, the $E$-connection is given by
\begin{eqnarray}
{}^E \nabla_{e} \alpha &:=& \calL_{\rho(e)} \alpha 
+ \bracket{\rho(\nabla e)}{\alpha}.
\label{Econoneform}
\end{eqnarray}
The $E$-connection on a $k$-form
$\alpha \in \Omega^k(M)= \Gamma(\wedge^k T^*M)$,
${}^E \nabla: 
\Gamma(\wedge^k T^*M) \rightarrow \Gamma(E^* \otimes \wedge^k T^*M)$
is given by applying the Leibniz rule recursively,
\begin{eqnarray}
{}^E \nabla_{e} (\beta \gamma) &=& 
({}^E \nabla_{e} \beta) \gamma - \beta ({}^E \nabla_{e} \gamma),
\label{Econoneform2}
\end{eqnarray}
for $\beta \in \Omega^1(M)$ and $\gamma \in \Omega^{k-1}(M)$.
Throughout our paper, the canonical $E$-connection is used
as the $E$-connection unless it is written.
%

Given an $E$-connection ${}^E \nabla$, the Lie algebroid differential ${}^E \rd$ can be generalized to the operation on the space $\Gamma(\wedge^k E' \otimes \wedge^m E^*)$, ${}^E \rd^{\nabla}: 
\Gamma(\wedge^k E' \otimes \wedge^m E^*)
\rightarrow \Gamma(\wedge^k E' \otimes \wedge^{m+1} E^*)$. For instance, see \cite{AbadCrainic}.
We set $E' = T^*M$ in the general formula
since only this case appears in our paper.
The concrete equation is as follows.
\begin{definition}
For $\Omega^k(M, \wedge^m E^*) = \Gamma(\wedge^k T^* M \otimes \wedge^m E^*)$,
the \textit{$E$-exterior covariant derivative}
${}^E \rd^{\nabla}: \Omega^k(M, \wedge^m E^*) \rightarrow \Omega^k(M, \wedge^{m+1} E^*)$ is defined by
\begin{eqnarray}
{}^E \rd^{\nabla} \alpha(e_1, \ldots, e_{m+1})
&:=& \sum_{i=1}^{m+1} (-1)^{i-1} 
{}^E \nabla_{e_i} 
(\alpha(e_1, \ldots, 
\check{e_i}, \ldots, e_{m+1}))
\nonumber \\ && 
+ \sum_{1 \leq i < j \leq m+1} (-1)^{i+j} \alpha([e_i, e_j], e_1, \ldots, \check{e_i}, \ldots, \check{e_j}, \ldots, e_{m+1}),
\label{Ediffconnection}
\end{eqnarray}
for $\alpha \in \Omega^k(M, \wedge^m E^*)$ and
$e_i \in \Gamma(E)$.
\end{definition}
Note that the $E$-exterior covariant derivative increases the order of $\wedge^m E^*$.

%

{
In Section \ref{sec:GNLSM}, we use the following $E$-covariant derivative,
\begin{eqnarray}
{}^E \nabla_{e} \alpha(e_1, \ldots, e_{m})
&=& (\calL^{\nabla}_{\rho(e)} \alpha) (e_1, \ldots, e_{m}),
\label{covariantLiederi}
\end{eqnarray}
where
\begin{eqnarray}
(\calL^{\nabla}_{\rho(e)} \alpha) 
=
(\iota_{\rho(e)} \nabla + \nabla \iota_{\rho(e)}) \alpha.
\end{eqnarray}
is the covariantized Lie derivative.
Comparing to Eq.~\eqref{LAdifferential},
$\calL^{\nabla}$  is also regarded as 
an extension of the anchor map $\rho$. We also use the following notation for
the covariantized exterior anchor map as
$\rho^{\nabla}: \Gamma(E) \times \Omega^k(M, \wedge^m E^*) 
 \rightarrow \Omega^k(M, \wedge^m E^*)$, 
$\rho^{\nabla}:(e, \alpha) \mapsto \rho^{\nabla}(e) \alpha$,
\begin{eqnarray}
&& (\rho^{\nabla}(e) \alpha) (e_1, \ldots, e_m) 
= (\calL^{\nabla}_{\rho(e)} \alpha)(e_1, \ldots, e_m),
\label{cov1}
 \end{eqnarray}
for $e \in \Gamma(E)$ and
$\alpha \in \Omega^k(M, \wedge^m E^*)$.
Note that $\rho^{\nabla}(e) f = \rho(e) f$ for $f \in C^{\infty}(M)$.
Moreover, it is useful to notice the manifestly covariant
expression of the Lie bracket on $\Gamma(E)$,
\begin{eqnarray}
&& [e_1, e_2] = -T(e_1, e_2)
+ \rho^{\nabla}(e_1) e_2 - \rho^{\nabla}(e_2) e_1.
\label{cov2}
\end{eqnarray}
where $T \in \Gamma(E \oplus \wedge^2 E^*)$ is the $E$-torsion defined by
\begin{eqnarray}
T(e_1, e_2) &:=& {}^E\nabla_{e_1} e_2 - {}^E\nabla_{e_2} e_1
- [e_1, e_2].
\label{Etorsion1}
 \end{eqnarray}
}
%
%
Eq.~\eqref{cov1} and the Lie bracket $[-,-]$
satisfy the following covariantized identities of the Lie algebroid,
\begin{eqnarray}
&& [\rho^{\nabla}(e_1), \rho^{\nabla}(e_2)] = \rho^{\nabla}([e_1, e_2]),
\label{covidentity1}
\\
&& [e_1, f e_2] = f [e_1, e_2] + {\rho^{\nabla}}(e_1)f \cdot e_2.
\label{covidentity2}
\end{eqnarray}
Eq.~\eqref{covidentity2} is easily obtained from 
${\rho^{\nabla}}(e_1)f = {\rho}(e_1)f$
and \eqref{covidentity1} is proved from \eqref{covidentity2}
and the Jacobi identity of the Lie bracket.

A boundary operator
on the exterior algebra of the vector bundle $E$, $\Gamma(\wedge^{\bullet} E)$,
is introduced.
Let $e_1 \wedge \ldots \wedge e_{m} \in \Gamma(\wedge^m E)$ be a decomposable element of $\Gamma(\wedge^m E)$ with $e_i \in \Gamma(E)$.
Then, a Lie algebroid homology operator 
is given by
\begin{align}
{}^E \partial_m (e_1 \wedge \ldots \wedge e_m)
&= \sum_{1 \leq i < j \leq m} (-1)^{i+j} [e_i, e_j] \wedge e_1 \wedge \ldots \wedge \check{e_i} \wedge \ldots \wedge \check{e_j} \wedge \ldots \wedge e_m,
\label{defhomology}
\end{align}
for $m > 1$.
Here $\Gamma(\wedge^0 E) \simeq C^{\infty}(M)$.
The corresponding homology operator on $\Gamma(\wedge^{\bullet} E)$ 
is defined by ${}^E \partial = \sum_{k=1}^{\mathrm{rk} E} \partial_m$.
The Lie algebroid homology operator ${}^E \partial$ 
satisfies $({}^E \partial)^2=0$.

One can refer to \cite{AbadCrainic, CrainicFernandes, DufourZung}
about the general theory of connections on a Lie algebroid.
Local coordinate expressions of formulas in a Lie algebroid 
are summarized in the Appendix.

\section{Homotopy momentum section}\label{sec:HMS}
In this section, we introduce a homotopy momentum section on a 
pre-multisymplectic manifold.

\begin{definition}
A pre-$n$-plectic manifold is a pair $(M, \momega)$, where $M$ is 
a smooth manifold and $\momega$ is a closed $(n+1)$-form on $M$.
\end{definition}
A pre-$n$-plectic manifold is also called a pre-multisymplectic manifold for $n \geq 2$.
If $\momega$ is nondegenerate, i.e., if $\iota_{v} \momega =0$ for a vector field $v \in \mathfrak{X}(M)$ implies $v=0$,
a pre-$n$-plectic manifold is called an $n$-plectic manifold.
A pre-$1$-plectic manifold is a pre-symplectic manifold
and a $1$-plectic manifold is nothing but a symplectic manifold.

We introduce a series of $k$-forms taking values in $\wedge^{n-k} E^*$, 
$\mu_k \in \Omega^k(M, \wedge^{n-k} E^*)
$, where $k=0,\ldots, n-1$.
Recall two derivations introduced in Section \ref{LAdiffandconn}.
The exterior covariant derivative  $\nabla$ 
and the $E$-exterior covariant derivative ${}^E \rd^{\nabla}$
act on the space $\Omega^k(M, \wedge^{n-k} E^*)$.
\begin{definition}\label{defhms}
A formal sum $\mu = \sum_{k=0}^{n-1} \mu_k$ is called a \textit{homotopy momentum section} if $\mu$ satisfies 
\begin{eqnarray}
&& (\nabla + {}^E \rd^{\nabla}) \mu = - \sum_{k=0}^{n} \iota_{\rho}^{n+1-k} \momega.
\label{homotopyMS}
\end{eqnarray}
\end{definition}
Here, ${\iota_{\rho}^k \omega} \in \Omega^{n-k}(M, \wedge^k E^*)$ is given by
\begin{eqnarray}
\iota_{\rho}^{k} \momega(v_{k+1}, \ldots, v_{n+1}) (e_1, \ldots, e_k)
&=& 
\iota_{\rho(e_1)} \ldots \iota_{\rho(e_k)} \momega(v_{k+1}, \ldots, v_{n+1}) 
\nonumber \\
&:=& 
\momega(\rho(e_k), \ldots, \rho(e_1), v_{k+1}, \ldots, v_{n+1}),
\end{eqnarray}
for $e_1, \ldots, e_k \in \Gamma(E)$ and 
$v_{k+1}, \ldots, v_{n+1} \in \mathfrak{X}(M)$.

Expanding the equation by form degree of both sides,
Equation \eqref{homotopyMS} is the following $n+1$ equations,
\begin{eqnarray}
&& \nabla \mu_{n-1} = - \iota_{\rho}^1 \momega,
\label{homotopyMS1}
\\
&& \nabla \mu_{n-2} + {}^E \rd^{\nabla} \mu_{n-1}= - \iota_{\rho}^{2} \momega,
\label{homotopyMS2}
\\
&& \vdots
\nonumber \\
&& \nabla \mu_{k-1} + {}^E \rd^{\nabla} \mu_{k}= - \iota_{\rho}^{n+1-k} \momega,
\label{homotopyMS3}
\\
&& \vdots
\nonumber \\
&& {}^E \rd^{\nabla} \mu_{0} = {}^E \rd \mu_{0} = - \iota_{\rho}^{n+1} \momega.
\label{homotopyMS4}
\end{eqnarray}
In general, neither the square of the operator $\nabla + {}^E \rd^{\nabla}$, $\nabla$, nor ${}^E \rd^{\nabla}$ is zero; i.e., they are not coboundary operators.
In this paper, we do not impose the condition such that
$\nabla^2 \mu = 0$ 
for the definition of a homotopy momentum section,
since these conditions are not needed for existence of a homotopy momentum section $\mu$. 

A homotopy Hamiltonian Lie algebroid is defined as follows.
\begin{definition}\label{hHamiltonianLA}
A Lie algebroid $E$ with a connection $\nabla$ and 
a homotopy momentum section 
$\mu \in \bigoplus_{k=0}^{n-1} \Omega^k(M, \wedge^{n-k} E^*)$
is called \textit{homotopy Hamiltonian} if 
Equation \eqref{homotopyMS} and $\nabla \iota_{\rho} \momega := {\nabla \iota_{\rho}^1 \momega} =0$ are satisfied.
\end{definition}

\begin{remark}
If the connection $\nabla$ is flat, i.e., $\nabla^2= 0$,
the $E$-covariant derivative satisfies $({}^E \rd^{\nabla})^2 =0$,
and $[\nabla, {}^E \rd^{\nabla} ]=0$,
$(\Omega^{\bullet}(M, \wedge^{\bullet} E^*), \nabla, {}^E \rd^{\nabla})$ gives
a double complex. However in general, both $\nabla$ and ${}^E \rd^{\nabla}$ are not necessarily differentials.

\end{remark}

\subsection{Examples}
In this section, some examples are listed.
Many geometric structures are regarded as an example of the homotopy momentum section.

\begin{example}[Momentum map]\label{mm}
A homotopy momentum section is a generalization of 
the momentum map on a symplectic manifold 
with a Lie group action. 

Let $(M, \momega)$ be a symplectic manifold with an action of a Lie group $G$ on $M$.
$\momega$ is a symplectic form, i.e., a nondegenerate closed $2$-form.
As explained in Example \ref{actionLA}, the Lie group action induces 
an action Lie algebroid structure on a trivial bundle 
$E = M \times \mathfrak{g}$ with a Lie algebra $\mathfrak{g}$ of $G$.
The anchor map is induced from the Lie algebra action 
$\mathfrak{g} \curvearrowright TM$.
Since $E = M \times \mathfrak{g}$ is the trivial bundle, we can take the trivial connection $\nabla = \rd$ with the de Rham differential $\rd$.
For $\alpha \in \Gamma(\wedge^m E^*) = C^{\infty}(M, \wedge^m \mathfrak{g}^*)$ and $e_i \in \mathfrak{g}$,
the $E$-exterior covariant derivative becomes
\begin{eqnarray}
{}^E \rd^{\nabla} \alpha(e_1, \ldots, e_{m+1}) 
&=& 
\mathrm{ad}^{*}_{\rho} \alpha(e_1, \ldots, e_{m+1}) 
+ \rd_{CE} \alpha(e_1, \ldots, e_{m+1}),
\label{adCE}
\end{eqnarray}
where $\mathrm{ad}^*$ is the action on $\alpha$ induced from the Lie algebra
action on $M$,
\begin{eqnarray}
\mathrm{ad}^*_{\rho} \alpha(e_1, \ldots, e_{m+1}) 
&:=& \sum_{i=1}^{m+1} (-1)^{i-1} \mathrm{ad}^*_{e_i} \alpha(e_1, \ldots, 
\check{e_i}, \ldots, e_{m+1})
\nonumber \\
&=& \sum_{i=1}^{m+1} (-1)^{i-1} \rho(e_i) \alpha(e_1, \ldots, 
\check{e_i}, \ldots, e_{m+1}),
\end{eqnarray}
for $e_i \in \mathfrak{g}$,
and $\rd_{CE}$ is the Chevalley-Eilenberg differential,
\begin{eqnarray}
\rd_{CE} \alpha(e_1, \ldots, e_{m+1}) 
&=& \sum_{1 \leq i < j \leq m+1} (-1)^{i+j} \alpha([e_i, e_j], e_1, \ldots, \check{e_i}, \ldots, \check{e_j}, \ldots, e_{m+1}).
\label{CEdifferential}
\end{eqnarray}
satisfying $(\rd_{CE})^2=0$.

Take $n=1$ in the definition of the homotopy momentum section.
A $1$-plectic form $\omega$ is nothing but a symplectic form,
and $\mu$ has only one component $\mu=\mu_0 \in C^{\infty}(M, \mathfrak{g}^*)$ 
which is a function taking a value in $\mathfrak{g}^*$.
Then, Equation \eqref{homotopyMS} reduces to two equations,
\begin{eqnarray}
&& \rd \mu_0 = - \iota_{\rho} \momega,
\label{MM1}
\\
&& {}^E \rd \mu_0 = - \iota_{\rho}^{ 2} \momega.
\label{MM2}
\end{eqnarray}
Substituting $e_1, e_2 \in \mathfrak{g}$ in Equation \eqref{MM2},
the left hand side is
\begin{eqnarray}
&& {}^E \rd \mu_0(e_1, e_2)
= \rho(e_1) \bracket{\mu_0}{e_2}
- \rho(e_2) \bracket{\mu_0}{e_1}
- \bracket{\mu_0}{[e_1, e_2]}.
\label{lefthand}
\end{eqnarray}
Using Equation \eqref{MM1} and 
$\rd \mu_0(e) = - \iota_{\rho} \momega(e) 
$, 
we obtain $\iota_{\rho(e_2)} \rd \mu_0(e_1)
= - \iota_{\rho(e_2)} \momega(\rho(e_1))$.
Using this equation, the second term of the right hand side 
in Eq.~\eqref{lefthand} becomes
\begin{eqnarray}
\rho(e_2) \bracket{\mu_0}{e_1}
&=& \iota_{\rho(e_2)} \rd \mu_0(e_1)
= - \iota_{\rho(e_2)} \momega(\rho(e_1))
= (\iota_{\rho}^{2} \momega)(e_1, e_2).
\label{righthand}
\end{eqnarray}
Therefore substituting \eqref{lefthand} and \eqref{righthand} to
\eqref{MM2}, we obtain the equation,
$\bracket{\mu_0}{[e_1, e_2]} = 
\rho(e_1) \bracket{\mu_0}{e_2}$.
i.e.,
\begin{eqnarray}
&& 
\mu_0([e_1, e_2]) = \mathrm{ad}^*_{e_1} \mu_0(e_2).
\label{MM3}
\end{eqnarray}
for $e_1, e_2 \in \mathfrak{g}$.
Therefore, Equation \eqref{homotopyMS} is equivalent to Equations
 \eqref{MM1} and \eqref{MM3}. 
Equation \eqref{MM3} means that the Lie algebra action is equivariant 
on $\mu_0$. Two equations shows that 
$\mu_0$ is a momentum map of the Lie group action on $M$.

Equation \eqref{MM3} can also be written as
\begin{eqnarray}
&& 
\rd_{CE} \mu_0(e_1, e_2) 
= - (\iota_{\rho}^{ 2} \momega)(e_1, e_2),
\label{MM4}
\end{eqnarray}
in view of Equation \eqref{MM1}.
\end{example}

\begin{example}[Momentum map on multisymplectic manifold]\label{mmmulti}
Let $n \geq 2$ and $(M, \momega)$ be an $n$-plectic manifold with 
an $n$-plectic form $\momega$.
Suppose an action of a Lie group $G$ on $M$ like the previous example \ref{mm}.
Similar to Example \ref{actionLA}, the Lie group action induces an action Lie algebroid on a trivial bundle $E = M \times \mathfrak{g}$, 
where $\mathfrak{g}$ is a Lie algebra of $G$.

If we set $\mu_k =0$ for $k=0,\ldots, n-2$, 
$\mu$ has only one component, an $n-1$ form $\mu=\mu_{n-1} \in 
\Omega^{n-1}(M, \mathfrak{g}^*)$. 
Since $\nabla=\rd$ again, 
Equation \eqref{homotopyMS} reduces to two equations,
\begin{eqnarray}
&& \rd \mu_{n-1} = - \iota_{\rho} \momega.
\label{MM21}
\\
&& {}^E \rd \mu_{n-1} = - \iota_{\rho}^{ 2} \momega.
\label{MM22}
\end{eqnarray}
Similar to Example \eqref{mm},
Equations \eqref{MM21} and \eqref{MM22} give the equation,
\begin{eqnarray}
&& 
\mu_{n-1}([e_1, e_2]) = \mathrm{ad}^*_{e_1} \mu_{n-1}(e_2).
\label{MM23}
\end{eqnarray}
for $e_1, e_2 \in \mathfrak{g}$.
Equations \eqref{MM21} and \eqref{MM22} are equivalent to Equations
 \eqref{MM21} and \eqref{MM23}. Two equations means that 
$\mu_{n-1}$ is a momentum map of the Lie group action on an $n$-plectic manifold $M$ \cite{Gotay:1997eg}.
\end{example}

\begin{example}[Momentum section]\label{ExMS}
We take $n=1$ and a general Lie algebroid $E$.
Since $n=1$, the base manifold $M$ is a pre-symplectic manifold.
Similar to the momentum map, a homotopy momentum section $\mu$ 
is only one component 
$\mu = \mu_0 \in \Gamma(E^*)$.
Equation \eqref{homotopyMS} reduces to two equations,
\begin{eqnarray}
&& \nabla \mu_0 = - \iota_{\rho} \momega,
\label{MS1}
\\
&& {}^E \rd \mu_0 = - \iota_{\rho}^{ 2} \momega.
\label{MS2}
\end{eqnarray}
Equations \eqref{MS1} and \eqref{MS2} are equivalent to 
the definition of a bracket-compatible momentum section 
on a Lie algebroid over a pre-symplectic manifold 
introduced by Blohmann and Weinstein \cite{Blohmann:2018}.
See also \cite{Kotov:2016lpx}.
Many examples of momentum sections which are not momentum maps have been discussed in \cite{Blohmann:2018}. 
One can refer to more examples of momentum sections appearing in physical theories in \cite{Ikeda:2019pef}.
\end{example}

In some cases, only parts of Equations \eqref{homotopyMS} in the definition of homotopy momentum sections appear.
Especially the $0$-form part of Equation \eqref{homotopyMS4} is important 
itself.

\begin{example}[Twisted Poisson structure]
We consider a twisted Poisson manifold $M$ in Example \ref{tPoisson}.
Then, as explained in Example \ref{tPoisson}, the twisted Poisson structure gives a Lie algebroid structure on $T^*M$.
Using the Lie algebroid differential ${}^E \rd$ given by this 
Lie algebroid, Equation \eqref{tPoisson1} of the twisted Poisson structure is rewritten as
\begin{eqnarray}
&& {}^E \rd \pi = - \iota_{\pi}^{ 3} H.
\end{eqnarray}
It means that $\mu_0 = \pi$ satisfies Equation \eqref{homotopyMS4}
with $H = \momega$ for $n=2$. Here $\mu_1=0$.
\end{example}

\begin{example}[twisted $R$-Poisson structure]
A twisted $R$-Poisson manifold $(M, \pi, J, H)$ \cite{Chatzistavrakidis:2021nom} is regarded as a higher dimensional generalization of the twisted Poisson structure.
Let $\pi \in \Gamma(\wedge^2 TM)$ be a Poisson bivector field on $M$. i.e.,
$\pi$ is a bivector field satisfying $[\pi,\pi]_S=0$.
$H$ is a closed $(n+1)$-form, and $J \in \Gamma(\wedge^{n} TM)$ is 
an $n$-multivector field.
\footnote{An $n$-multivector field is denoted by $J$ in this paper though
it is denoted by $R$ in \cite{Chatzistavrakidis:2021nom}.
}
The $R$-Poisson structure $(\pi, J, H)$ is defined by 
the triple satisfying the equation
\begin{eqnarray}
&& [\pi, J]_S = \frac{}{} 
\bracket{\otimes^{n+1} \pi}{H}.
\label{RPoisson2}
\end{eqnarray}
The Poisson bivector field $\pi$ induces a Lie algebroid structure
on $T^*M$ as written in Example \ref{Poisson}.
Equation \eqref{RPoisson2} is equivalent to the equation,
\begin{eqnarray}
&& {}^E \rd J = - \iota_{\pi}^{n+1} H,
\label{Rtwisted}
\end{eqnarray}
where the Lie algebroid differential ${}^E \rd$ is one
with respect to the Lie algebroid $T^*M$ induced from the Poisson structure.
The condition \eqref{Rtwisted} is nothing but Equation \eqref{homotopyMS4}
with $\omega = H$, $\mu_0 = J$ and $\mu_k=0$ for $k=1,\ldots,n-1$.
The twisted $R$-Poisson structure is also regarded as a special case of the homotopy momentum section.
Refer to \cite{Ikeda:2021rir}
for a generalization of the twisted $R$-Poisson structure 
to a general Lie algebroid and a corresponding topological sigma model.
\end{example}

\section{Weak homotopy momentum section
and equivariant $E$-differential form}\label{sec:WHMS}
In this section, we analyze the Lie algebroid over an $n$-plectic manifold in details for $n \geq 2$, and consider a generalization of the weak homotopy momentum map \cite{Herman:2018box}. Moreover we introduce a condition, \textit{equivariant}, which is nontrivial only for $n \geq 2$.

A pairing of $\Gamma(\wedge^{\bullet} E)$ and
$\Gamma(\wedge^{\bullet} E^*)$, 
$\bracket{-}{-}: \Gamma(\wedge^{\bullet} {E}) \times \Gamma(\wedge^{\bullet} E^*) \rightarrow \bR$ is induced from the natural pairing of
$E$ and $E^*$. 
Let $u_i$ and $u^i$ be basis on $E$ and $E^*$.
Then,
for an $E$-differential form $\alpha = \frac{1}{m!} \alpha_{a_1 \dots a_m} 
{u^{a_1} \wedge \ldots \wedge u^{a_m}} \in \Gamma(\wedge^m E^*)$
and a decomposable element of $\Gamma(\wedge^{\bullet} E)$,
$e_1 \wedge \ldots \wedge e_m 
= \frac{1}{m!} e_1^{a_1} \ldots e_m^{a_m} u_{a_1} 
\wedge \ldots \wedge u_{a_m}$ for $e_i^{a_i} \in \bR$, 
the pairing is given by
\begin{align}
\alpha(e_1, \ldots, e_m) &= 
\bracket{\alpha}{e_1, \ldots, e_m} 
= \frac{1}{m!}  \alpha_{a_1 \dots a_m} e_1^{a_1} \ldots e_m^{a_m}.
\end{align}

An $m$-th \textit{Lie kernel} $\calP_m \subset \Gamma(\wedge^{m} E)$ 
for a Lie algebroid is a subset of $\Gamma(\wedge^{m} E)$
defined by the kernel of ${}^E \partial_m$, 
\begin{align}
\calP_m &= \{w \in \Gamma(\wedge^{m} E)
|w \in \mathrm{Ker}({}^E \partial_m) \},
\end{align}
where ${}^E \partial_m$ is the Lie algebroid homology operator defined in
Equation \eqref{defhomology}.
A (total) Lie kernel is defined by $\calP = \oplus_{k=1}^{\mathrm{rank} E} \calP_k$.

Let $\mu = \sum_{k=0}^{n-1} \mu_k$, where $\mu_k$ be a map $\mu_k: \calP_{n-k} \rightarrow \Omega^k(M)$.
{An element $\mu \in \Omega^k(M, \wedge^{n-k} E^*)$ is also regarded 
as the map $\mu_k: \calP_{n-k} \rightarrow \Omega^k(M)$.}
We define a generalization of a weak homotopy momentum map.
\begin{definition}\label{defwhMS}
Let $E$ be a Lie algebroid over a pre-$n$-plectic manifold $M$.
$\mu: \calP \rightarrow \sum_{k=0}^{n-1} \Omega^k(M)$ 
is called a \textit{weak homotopy momentum section} if $\mu$ satisfies
\begin{eqnarray}
&& \nabla \mu = - \sum_{k=1}^{n-1} \iota_{\rho}^{n+1-k} \momega.
\label{whomotopyMS}
\end{eqnarray}
\end{definition}

Next, we introduce the notion of equivariance for sections of $\wedge^m E^*$ as follows:
\begin{definition}\label{defequivariant}
$\alpha \in \Omega^l(M, \wedge^m E^*)$ is called \textit{equivariant} 
with respect to a Lie algebroid $(E, \rho, [-,-])$
if it satisfies 
\begin{align}
{}^E \nabla_e (\alpha(e_1, \ldots, e_m))
&= \sum_{i=1}^m (-1)^{i-1} \alpha({[e, e_i]}, e_1, \ldots, \check{e_i}, \ldots, e_m),
\label{equivequ}
\end{align}
for $e, e_i \in \Gamma(E)$
\footnote{For $m=0$, the righthand side of Eq.~\eqref{equivequ} is understood as zero.}.
\end{definition}

In the action Lie algebroid, the connection is trivial $\nabla = \rd$.
Thus Equation \eqref{equivequ} reduces to
\begin{align}
\rho(e) \alpha(e_1, \ldots, e_m)
&= \sum_{i=1}^m (-1)^{i-1} \alpha([e, e_i], e_1, \ldots, \check{e_i}, \ldots, e_m).
\label{equivequ2}
\end{align}
For a momentum map $\mu$,
Equation \eqref{equivequ2} follows from Equation \eqref{MM3} in Example \ref{mm}.
\begin{proposition}
A momentum map $\mu$ on a symplectic manifold $M$ with a Lie group action is equivariant in the sense of Definition \ref{defequivariant}.
\end{proposition}
Let $M$ be a general $n$-plectic manifold and $E$ be a Lie algebroid.
Then if $\alpha$ is equivariant, applying Eq.~\eqref{equivequ},
the $E$-exterior covariant derivative \eqref{LAdifferential} is simplified to
\begin{align}
{}^E \rd^{\nabla} \alpha(e_1, \ldots, e_{m+1}) 
&= - \sum_{1 \leq i < j \leq m+1} (-1)^{i+j} \alpha({[e_i, e_j]}, e_1, \ldots, \check{e_i}, \ldots, \check{e_j}, \ldots, e_{m+1}).
\label{equivariantEq}
\end{align}
Thus, the following formula holds.
\begin{lemma}\label{equivLK}
Let $\alpha \in \Gamma(\wedge^m E^*)$ be equivariant.
Then, for a decomposable element of a Lie kernel,
$e_1 \wedge \ldots \wedge e_{m+1} \in \calP_{m+1}$,
\begin{eqnarray}
{}^E \rd^{\nabla} \alpha (e_1, \ldots, e_{m+1}) =0.
\label{kerEd}
\end{eqnarray}
\end{lemma}
In fact, the left hand side of \eqref{kerEd}
is Equation \eqref{equivariantEq}, and 
it is zero since 
$e_1 \wedge \ldots \wedge e_m$ is $\partial$-closed since it is a decomposable element of the Lie kernel.
The lemma holds since general elements of $\calP_{m+1}$ are linear combinations of decomposable elements.

Compare to the definition of the homotopy momentum section 
\eqref{homotopyMS} and the weak homotopy momentum section \eqref{whomotopyMS}.
A homotopy momentum map and a weak homotopy momentum map
 have been compared in \cite{Mammadova-Ryvkin}.
If $\mu \in \bigoplus_{k=0}^{n-1} \Omega^k(M, \wedge^{n-k} E^*)$ satisfies Equation \eqref{homotopyMS} and is equivariant,
$\mu$ satisfies Equation \eqref{whomotopyMS} from Lemma \ref{equivLK}.
i.e.,
\begin{proposition}
Let a homotopy momentum section $\mu$ be equivariant, i.e., 
satisfies Eq.~\eqref{equivequ}.
Then, $\mu$ is a weak homotopy momentum section.
\end{proposition}

\section{Relation to homotopy momentum map and weak homotopy momentum map}\label{sec:HMM}
In this section, we discuss the relation between 
a homotopy momentum section in this paper and a homotopy momentum map in \cite{Fregier:2013dda}.
The correspondence of weak versions is also discussed.

Let $(M, \momega)$ be an $n$-plectic manifold.
Suppose an action of a Lie group $G$ on $M$. 
The action of $G$ induces the corresponding infinitesimal Lie algebra action on $M$ as vector fields, $\rho:\mathfrak{g} \rightarrow TM$.
We consider two differentials on the space,
$\Omega^{k}(M, \wedge^{n-k} \mathfrak{g}^*) =
\Omega^{k}(M) \otimes \wedge^{n-k} \mathfrak{g}^*$,
$\rd := \rd \otimes 1$ and $\rd_{CE} := 1 \otimes \rd_{CE}$
where $\rd$ is the de Rham differential on $\Omega^{k}(M)$ and 
$\rd_{CE}$ is the Chevalley-Eilenberg differential on $\wedge^{n-k} \mathfrak{g}^*$.
Introduce a $k$-form taking a value in $\wedge^{n-k} \mathfrak{g}^*$,
$\hmu_k \in \Omega^{k}(M, \wedge^{n-k} \mathfrak{g}^*)$, 
where $k=0, \ldots, n-1$,
and their formal sum, $\hmu = \sum_{k=0}^{n-1} \hmu_k$.
\begin{definition}
A homotopy momentum map $\hmu$ is defined by
\footnote{Sign factors of the equation are different from \cite{Fregier:2013dda}. Sign factors are arranged to original ones by proper redefinitions.}
\begin{eqnarray}
&& (\rd + \rd_{CE}) \hmu = \sum_{k=0}^{n} (-1)^{n-k+1} \iota_{\rho}^{n+1-k} \momega.
\label{homotopyMM}
\end{eqnarray}
\end{definition}
Expanded by the form degree, 
Equation \eqref{homotopyMM} becomes the following $n$ equations,
\begin{eqnarray}
&& \rd \hmu_{n-1} = - \iota_{\rho}^1 \momega,
\label{homotopyMM1}
\\
&& \rd \hmu_{n-2} + \rd_{CE} \hmu_{n-1}= \iota_{\rho}^{2} \momega,
\label{homotopyMM2}
\\
&& \vdots
\nonumber 
\\
&& \rd \hmu_{k-1} + \rd_{CE} \hmu_{k}= (-1)^{n-k+1} \iota_{\rho}^{n+1-k} \momega,
\label{homotopyMM3}
\\
&& \vdots
\nonumber 
\\
&& \rd_{CE} \hmu_{0} = (-1)^{n+1} \iota_{\rho}^{n+1} \momega.
\label{homotopyMM4}
\end{eqnarray}
For $n=1$, a homotopy momentum map is equivalent to a momentum map.
In fact, in this case, $\momega$ is a symplectic form and
$\hmu = \hmu_0 \in C^{\infty}(M, \mathfrak{g}^*)$ has only one element. 
Equation \eqref{homotopyMM} reduces to two equations,
\begin{eqnarray}
\rd \hmu_{0} &=& - \iota_{\rho} \momega,
\label{homotopyMM21} \\
\rd_{CE} \hmu_{0} &=& \iota_{\rho}^{2} \momega.
\label{homotopyMM24}
\end{eqnarray}
Let $X_e = \rho(e)$ be a vector field on $M$ given by an element $e$ of the Lie algebra $\mathfrak{g}$.
Since $\rd_{CE} \hmu_{0}(e_1, e_2) = - \hmu_0([e_1, e_2])$
and 
$\iota_{\rho}^{2} \momega(e_1, e_2) = 
- \momega(\rho(e_1), \rho(e_2)) = - \momega(X_1, X_2)$,
Equation \eqref{homotopyMM24} is equivalent to
\begin{eqnarray}
\hmu_{0}([e_1, e_2]) &=& \momega(X_1, X_2),
\label{homotopyMM34}
\end{eqnarray}
which means that $\hmu_{0}$ is an (anti)-homomorphism from 
$\mathfrak{g}$ to $TM$.
Equations \eqref{homotopyMM21} and \eqref{homotopyMM34} show that $\hmu_0$ is a momentum map on the symplectic manifold $M$.
For $n=1$,
both a homotopy momentum map and a homotopy momentum section
in Example \ref{mm} (a momentum section) give a momentum map on a symplectic manifold.
Therefore, the result is summarized as follows.
\begin{proposition}
For $n=1$, $M$ is a symplectic manifold. 
Then, a homotopy momentum section on an action Lie algebroid $M \times \mathfrak{g}$ is equivalent to a homotopy momentum map and a momentum map on $M$
with the Lie algebra action $\mathfrak{g}$.
\end{proposition}

We consider the $n \geq 2$ case. 
As we see in Example \ref{mm}, the Lie group action induces 
an action Lie algebroid structure on the trivial bundle
$E = M \times \mathfrak{g}$.
Then, 
the $E$-exterior covariant derivative becomes
${}^E \rd^{\nabla} = \mathrm{ad}^*_{\rho} + \rd_{CE}$
from Equation \eqref{adCE}.

Suppose that the homotopy momentum section $\mu$ on the action Lie algebroid is equivariant. Then, from Equation \eqref{equivariantEq}, we obtain ${}^E \rd^{\nabla} = - \rd_{CE}$.
Thus, Equation \eqref{homotopyMS} of the homotopy momentum section reduces to 
\begin{eqnarray}
&& (\rd - \rd_{CE}) \mu = - \sum_{k=0}^{n} \iota_{\rho}^{n+1-k} \momega.
\label{homotopyMSAL}
\end{eqnarray}
If we take $\hmu_k = (-1)^{n-k+1} \mu_k$, 
Equation \eqref{homotopyMSAL} coincides with
Equation \eqref{homotopyMM1}--\eqref{homotopyMM4}.
Therefore we obtain the following relation between a homotopy momentum section and a homotopy moment map for a Lie algebra action.
\begin{theorem}
Let $n \geq 2$ and a Lie algebroid be the action Lie algebroid $E = M \times \mathfrak{g}$ with a Lie algebra $\mathfrak{g}$. Then, 
If the homotopy momentum section $\mu$ for the action Lie algebroid $E$ 
is equivariant,
it is a homotopy momentum map for $\mathfrak{g}$.
\end{theorem}

Next, we discuss relations between a weak homotopy momentum section and 
a weak homotopy momentum map.
Let $\wedge^m \mathfrak{g}$ be an exterior algebra of the Lie algebra $\mathfrak{g}$.
A Lie algebra homology operator $\partial_m: \Gamma(\wedge^m \mathfrak{g}) \rightarrow \Gamma(\wedge^{m-1} \mathfrak{g})$ is defined by
\begin{align}
\partial_m (e_1 \wedge \ldots \wedge e_m)
&= \sum_{1 \leq i < j \leq m} (-1)^{i+j} ([e_i, e_j] \wedge e_1 \wedge \ldots \wedge \check{e_i} \wedge \ldots \wedge \check{e_j} \wedge \ldots \wedge e_m).
\label{Liealgebrahomology}
\end{align}
for a decomposable element $e_1 \wedge \ldots \wedge e_{m} \in \wedge^m \mathfrak{g}$ with $e_i \in \mathfrak{g}$.
The Lie algebra homology operator on $\wedge^{\bullet} \mathfrak{g}$ is defined by $\partial = \sum_m \partial_m$.

An $m$-th Lie kernel $\calP_m \subset \Gamma(\wedge^{m} \mathfrak{g})$ 
for a Lie algebra $\mathfrak{g}$
is defined by the kernel of a Lie algebra homology operator $\partial_m$,
\begin{align}
\calP_m &= \{w \in \Gamma(\wedge^{m} \mathfrak{g})
|w \in \mathrm{Ker}(\partial_m) \},
\end{align}
and a (total) Lie kernel is defined by $\calP = \oplus_{m=1}^{\mathrm{dim} \mathfrak{g}} \calP_m$.

Consider a map $\hmu_k: \calP \rightarrow \Omega^{k}(M)$,
where $k=0, \ldots, n-1$.
The map $\hmu_k: \calP \rightarrow \Omega^{k}(M)$ is regarded as $\hmu_k \in \Omega^{k}(M, \wedge^{n-k} \calP^*)$.
For the formal sum of $\hmu_k$, $\hmu = \sum_{k=0}^{n-1} \hmu_k$,
a weak homotopy momentum map $\hmu$ is defined by \cite{Herman:2017}
\begin{eqnarray}
&& \rd \hmu = \sum_{k=0}^{n} (-1)^{n-k+1} \iota_{\rho}^{n+1-k} \momega,
\label{whomotopyMM}
\end{eqnarray}
A homotopy momentum map $\hmu$ is always a weak homotopy momentum map
since $\rd_{CE} \hmu =0$ if the homotopy momentum map $\hmu$ is restricted to $\calP$.

We consider a weak homotopy momentum section in the definition \ref{defwhMS}
on the action Lie algebroid.
A homotopy momentum section $\mu$ is regarded as a map from the Lie kernel $\calP$ to $\Omega^{\bullet}(M)$. A connection is $\nabla = \rd$ on the action Lie algebroid.
Then, the Lie algebroid homology operator 
${}^E \partial_m$ in Equation \eqref{defhomology} with $\nabla = \rd$
is the same as the Lie algebra homology operator $\partial_m$
\eqref{Liealgebrahomology}.
Thus Equation \eqref{whomotopyMS} and \eqref{whomotopyMM} is equivalent 
if we et
$\hmu_k = (-1)^{n-k+1} \mu_k$.
Therefore, we obtain 
\begin{theorem}
Let $n \geq 2$ and a Lie algebroid be the action Lie algebroid. Then, if $\mu$ is a weak homotopy momentum section, it is a weak homotopy momentum map.
\end{theorem}

\section{Gauged nonlinear sigma model with Wess-Zumino term
}\label{sec:GNLSM}
\noindent
In this section, we show that a homotopy momentum section appears in 
the Lie algebroid gauging of an $n$-dimensional nonlinear sigma model with 
Wess-Zumino (WZ) term.
A nonlinear sigma model is a mechanics on the mapping space from an $(n+1)$-dimensional smooth manifold $\Xi$ with $n$-dimensional boundary $\Sigma = \partial \Xi$ to a $d$-dimensional manifold $M$.
Suppose that a target space $M$ is a $d$-dimensional pre-$n$-plectic 
Riemannian manifold with a metric $g$.
Moreover suppose a Lie algebroid $E$ over $M$.
The claim is that if gauge invariance of the action functional is required, we need an additional condition, existence of homotopy momentum sections on $E$.

Suppose a metric $g$, and a closed $(n+1)$-form $H$, which is a pre-$n$-plectic from.
Next we introduce an $(n+1)$-dimensional worldvolume $\Xi$ with boundary $\Sigma = \partial \Xi$, and a metric $\gamma$ on $\Sigma$.
Let $X:\Xi \rightarrow M$ be a map from $\Xi$ to $M$.
The action functional of the $n$ dimensional 
nonlinear sigma model with a Wess-Zumino term is as follows,
\begin{align}
S &= \int_{\Sigma} \frac{1}{2} X^* g(\rd X, *\rd X)
+ \int_{\Xi} X^* H
\nonumber \\ 
&= \int_{\Sigma} \frac{1}{2} g_{ij}(X) \rd X^i \wedge * \rd X^j
+ \int_{\Xi} \frac{1}{(n+1)!} H_{i_1\ldots i_{n+1}}(X) \rd X^{i_1} \wedge \ldots \wedge \rd X^{i_{n+1}}
\nonumber \\ 
&= \int_{\Sigma} d^n \sigma \frac{1}{2} g_{ij}(X) \sqrt{\gamma} 
\gamma^{\mu\nu} \partial_{\mu} X^i \partial_{\nu} X^j
\nonumber \\ & \qquad
+ \int_{\Xi} \frac{1}{(n+1)!} d^{n+1} \sigma H_{i_1\ldots i_{n+1}}(X) 
\epsilon^{M_1\ldots M_{n+1}} \partial_{M_1} X^{i_1} 
\wedge \ldots \wedge \partial_{M_{n+1}} X^{i_{n+1}},
\label{ndsigmamodel}
\end{align}
where $i, j$ are indices for local coordinates on $M$,
$\mu, \nu, \ldots$ are indices for local coordinates on $\Sigma$,
$M, N$ are indices for local coordinates on $\Xi$.
$\sigma^M$ is the local coordinate on $\Xi$.
$X^*g = g(X)$ is a pullback of the metric $g$ and
$X^* H = H(X) = \frac{1}{(n+1)!} H_{i_1\ldots i_{n+1}}(X) \rd X^{i_1} \wedge \ldots \wedge \rd X^{i_{n+1}}$ in the second term is a pullback of the pre-$n$-plectic form $H$ on $M$.
$*$ is the Hodge star with respect to the metric $\gamma$ on $\Sigma$. $n=2$ case is most important since it is the string sigma model with NS-flux $3$-form $H$.
Mainly we take expressions with local coordinate indices of the target space $M$ like the second expression in Equation \eqref{ndsigmamodel}.

If $M$ has an isometry with respect to a Lie group $G$ action, 
its infinitesimal Lie algebra $\mathfrak{g}$ acts on 
$M$ as a set of vector fields given by a homomorphism 
$\rho:M \times \mathfrak{g} \rightarrow TM$.
Then, the action of the Lie algebra on $X$ is given by $\rho$ as
\begin{eqnarray}
\delta X^i &=& \rho^i_a(X) \epsilon^a,
\label{transformationofX} 
\end{eqnarray}
where $\rho:= \rho^i_a(x) e^a \partial_i$ for
the basis $e^a$ of $\mathfrak{g}^*$.
$\epsilon^a$ is a gauge parameter, which is a function on $\Xi$ taking a value in $\mathfrak{g}$.
Since $\rho$ must be a Lie algebra homomorphism, it satisfies
Equation \eqref{almostLA} in the action Lie algebroid.

Gauge invariance of the action functional $S$ under 
the transformation \eqref{transformationofX} imposes
the following conditions of $H$ and $g$,
\begin{eqnarray}
\calL_{\rho(\epsilon)} g &=& 0,
\label{killingg2}
\\
\calL_{\rho(\epsilon)} H &=& \rd \beta(\epsilon),
\label{killlingB2}
\end{eqnarray}
where $\beta$ is an arbitrary $n$-form taking a value in $E^*$.

Next, we consider \textit{gauging} of the $n$-dimensional sigma model 
\eqref{ndsigmamodel}
by introducing a connection $D$ on $\Sigma$ and its connection
$1$-form $A$.
Gauging is a deformation of the theory by introducing 
connections (gauge fields) on the worldvolume.
The trivial bundle $E= M \times \mathfrak{g}$ is replaced to 
a general vector bundle $E$ with the fiber $\mathfrak{g}$.
A Lie algebra structure on $\mathfrak{g}$, which gives an action Lie algebroid on $E= M \times \mathfrak{g}$, is replaced to a Lie algebroid structure on $E$.

Introduce a $1$-form $A \in \Omega^1(\Sigma, X^*E)$ on $\Sigma$ taking a value in the pullback of $E$.
Moreover introduce a connection $\nabla$ on $E$ and its connection $1$-form 
$\omega_{ai}^b(x) \rd x^i \otimes e^a \otimes e_b$.
Take local coordinate expressions of the anchor map and the Lie bracket as
\begin{eqnarray}
\rho(e_a) &:=& \rho^i_a(x) \partial_i
\\
~[e_a, e_b]&:=& C_{ab}^c(x) e_c.
\end{eqnarray}
We consider the following target space covariant gauge transformations,
\footnote{We can consider more general gauge transformations, for instance, by adding the term $\phi_{bi}^a(X) \epsilon^b * D X^i$
to the right hand side of Equation \eqref{gaugetransformationA}
\cite{Chatzistavrakidis:2016jfz, Chatzistavrakidis:2016jci, Chatzistavrakidis:2017tpk}. We set this term zero in this paper.}
\begin{eqnarray}
\delta X^i &=& \rho^i_a(X) \epsilon^a,
\label{gaugetransformationX}
\\
\delta A^a &=& 
\rd \epsilon^a + C_{bc}^a(X) A^b \epsilon^c
+ \omega_{bi}^a(X) \epsilon^b D X^i.
\label{gaugetransformationA}
\end{eqnarray}
Here $D X^i = \rd X^i - \rho^i_a(X) A^a$ is the 
covariant derivative of $X^i$ on $\Sigma$. In fact, Equations
\eqref{gaugetransformationX} and \eqref{gaugetransformationA} are covariant under diffeomorphisms and coordinate transformations on the target bundle $E$
\cite{Mayer:2009wf}.

In order to make the action functional to be invariant under
gauge transformations \eqref{gaugetransformationX} 
and \eqref{gaugetransformationA}, we deform the action 
functional $S$ \eqref{ndsigmamodel} adding some terms
based on so called following Hull-Spence type ansatz 
\cite{Hull:1989jk, Hull:1990ms},
\begin{eqnarray}
S_{LA} = S_g + S_H + S_{\mu},
\label{ndgaugedsigmamodel}
\end{eqnarray}
where 
\begin{eqnarray}
S_g &=& \int_{\Sigma} \frac{1}{2} g_{ij} DX^i \wedge *DX^j
\\
S_H &=& \int_{\Xi} \frac{1}{(n+1)!} H_{i_1\ldots i_{n+1}}(X) 
\rd X^{i_1} \wedge \ldots \wedge \rd X^{i_{n+1}},
\\
S_{\mu} &=& \int_{\Sigma} 
\sum_{k=0}^{n-1} \frac{1}{k!(n-k)!} \tmu^{(k)}_{i_1 \ldots i_k a_{k+1} \ldots a_{n}}
(X) \rd X^{i_1} \wedge \ldots \wedge \rd X^{i_k} 
\wedge A^{a_{k+1}} \wedge \ldots \wedge A^{a_{n}},
\label{gaugedansatz}
\end{eqnarray}
where $\tmu^{(k)}$ is the pullback of a $k$-form on $M$ taking a value in 
$\wedge^{n-k}E^*$, $\tmu^{(k)} \in X^* \Omega^k(M, \wedge^{n-k}E^*)$.
In our case, the gauge structure is not a Lie algebra but a Lie algebroid.
Requiring \eqref{ndgaugedsigmamodel} is invariant under
gauge transformations \eqref{gaugetransformationX} 
and \eqref{gaugetransformationA}, we obtain geometric conditions for 
a metric $g$, $H$ and $\tmu^{(k)}$.
\footnote{For simplicity of notation, we use the same notation for geometric quantities on $M$ and $E$ and their pullbacks by the map $X$.}

The condition of the metric $g$
is independent of conditions for $H$ and $\tmu^{(k)}$
since $\delta S_g=0$ is required.
We obtain the following condition for $g$,
\begin{eqnarray}
&& {}^E\nabla g=0,
\label{conditionofg}
\end{eqnarray}
where ${}^E\nabla$ is the $E$-connection induced on $TM \otimes TM$.
We can check that if the connection is trivial, it Equation \eqref{conditionofg} reduces to Equation \eqref{killingg2}.

Next we fix conditions of $H$ and $\tmu^{(k)}$
to make the action functional \eqref{ndgaugedsigmamodel} gauge invariant.
This condition is $\delta (S_H+ S_{\mu})=0$.
The gauge transformation of $S_{\mu}$ is
\begin{align}
\delta S_{\mu} &= \int_{\Sigma} 
\sum_{k=0}^{n} \frac{1}{k!(n-k)!} 
\epsilon^b \rho^j_b \partial_j \tmu^{(k)}_{i_1 \ldots i_k a_{k+1} \ldots a_{n}}
(X) \rd X^{i_1} \wedge \ldots \wedge \rd X^{i_k} 
\wedge A^{a_{k+1}} \wedge \ldots \wedge A^{a_{n}}
\nonumber \\ & 
\quad + \sum_{k=1}^{n} \frac{1}{(k-1)!(n-k)!} 
\epsilon^b \partial_{i_1} \rho^j_b \tmu^{(k)}_{ji_2 \ldots i_k a_{k+1} \ldots a_{n}}
(X) \rd X^{i_1} \wedge \ldots \wedge \rd X^{i_k} 
\wedge A^{a_{k+1}} \wedge \ldots \wedge A^{a_{n}}
\nonumber \\ & 
\quad + \sum_{k=1}^{n} \frac{1}{(k-1)!(n-k)!} 
\rho^j_b \tmu^{(k)}_{ji_2 \ldots i_k a_{k+1} \ldots a_{n}}
(X) \rd \epsilon^b \wedge \rd X^{i_2} \wedge \ldots \wedge \rd X^{i_k} 
\wedge A^{a_{k+1}} \wedge \ldots \wedge A^{a_{n}}
\nonumber \\ & 
\quad + \sum_{k=0}^{n-1} (-1)^k \frac{1}{k!(n-k-1)!} 
\tmu^{(k)}_{i_1 \ldots i_k b a_{k+2} \ldots a_{n}}(X) 
\rd \epsilon^b \wedge \rd X^{i_1} \wedge \ldots \wedge \rd X^{i_k} 
\wedge A^{a_{k+2}} \wedge \ldots \wedge A^{a_{n}}
\nonumber \\ & 
\quad + \sum_{k=0}^{n-1} \frac{1}{k!(n-k-1)!} 
\tmu^{(k)}_{i_1 \ldots i_k c a_{k+2} \ldots a_{n}} C_{a_{k+1}b}^c
\epsilon^b \rd X^{i_1} \wedge \ldots \wedge \rd X^{i_k} 
\wedge A^{a_{k+1}} \wedge \ldots \wedge A^{a_{n}}
\nonumber \\ & 
\quad + \sum_{k=0}^{n-1} \frac{1}{k!(n-k-1)!} 
\tmu^{(k)}_{i_1 \ldots i_k c a_{k+2} \ldots a_{n}} \omega_{bi_{k+1}}^c
\epsilon^b \rd X^{i_1} \wedge \ldots \wedge \rd X^{i_{k+1}} 
\wedge A^{a_{k+2}} \wedge \ldots \wedge A^{a_{n}}
\nonumber \\ & 
\quad + \sum_{k=0}^{n-1} \frac{1}{k!(n-k-1)!} 
(-1) \tmu^{(k)}_{i_1 \ldots i_k c a_{k+2} \ldots a_{n}} \omega_{bj}^c
\rho^j_{a_{k+1}} \epsilon^b \rd X^{i_1} \wedge \ldots \wedge \rd X^{i_k} 
\wedge A^{a_{k+1}} \wedge \ldots \wedge A^{a_{n}}.
\label{gaugeSmu}
\end{align}
Since $\delta H = \calL_{\rho(\epsilon)} H$
from the gauge transformation of X, \eqref{gaugetransformationX},
the gauge transformation of $S_H$ is
\begin{eqnarray}
\delta S_H
&=& \int_{\Xi} \calL_{\rho(\epsilon)} H
=  \int_{\Xi} \rd \iota_{\rho(\epsilon)} H
= \int_{\Sigma} \iota_{\rho(\epsilon)} H
\nonumber \\ &=& 
\frac{1}{n!} \int_{\Sigma} \epsilon^a 
\rho^j_a H_{ji_2\ldots i_{n+1}}(X) \rd X^{i_2} \wedge \ldots 
\wedge \rd X^{i_{n+1}},
\label{gaugeSH}
\end{eqnarray}
where $\rd H=0$ and $\partial \Xi = \Sigma$ are used. 
Since the sum of coefficient functions of 
Equations \eqref{gaugeSmu} and \eqref{gaugeSH} must vanish,
we obtain the following conditions for $H$ and $\tmu^{(k)}$.
From terms including $\rd \epsilon^b$, we obtain
the algebraic condition for $\tmu^{(k-1)}$ and $\tmu^{(k)}$
for $k=1,\ldots, n-1$,
\begin{eqnarray}
&& \tmu^{(k-1)}(e_k, \ldots, e_n) = (-1)^k \iota_{\rho(e_k)} 
\tmu^{(k)}(e_{k+1}, \ldots, e_n) + \mathrm{Cycl}(e_k, \ldots, e_n).
\label{conditionofgsm21}
\end{eqnarray}
Comparing coefficients of terms,
$\epsilon^b \rd X^{i_1} \wedge \ldots \wedge \rd X^{i_k} 
\wedge A^{a_{k+1}} \wedge \ldots \wedge A^{a_{n}}$,
we obtain differential equations for $\tmu^{(k)}$ and $H$,
\begin{align}
&
\iota_{\rho(e)} H + \calL_{\rho(e)} \tmu^{(n)}
+ \bracket{\omega(e)}{\tmu^{(n-1)}} = 0,
\qquad (k=n)
\label{coef1}
\\
&
\sum_{l=k+1}^n (-1)^{l-k-1} \calL_{\rho(e_l)} \tmu^{(k)}(e_{k+1}, \ldots, \check{e}_l, \ldots, e_n)
\nonumber \\
&
- \sum_{k+1 \leq l < m \leq n} (-1)^{l+m-1}
\tmu^{(k)}([e_l, e_m], e_{k+1}, \ldots, \check{e}_l, \ldots, \check{e}_m, \ldots, e_n)
\nonumber \\ &
+ \sum_{l=k}^n (-1)^{l-k} \bracket{\omega(e_l)}{\tmu^{(k-1)}}
(e_{k}, \ldots, \check{e}_l, \ldots, e_n)
\nonumber \\ &
- \sum_{l \neq m} (-1)^{l+m-1} \bracket{\iota_{\rho(e_m)} \omega(e_l)}{\tmu^{(k)}}(e_{k+1}, \ldots, \check{e}_l, \ldots, \check{e}_m, \ldots, e_n)
= 0,
\qquad (k=1, \cdots, n-1)
\label{coef2}
\\
&
{}^E \rd \tmu^{(0)}(e_1, \ldots, e_n)
- \sum_{l=1}^n (-1)^{l-1} \bracket{\iota_{\rho(e_l)}\omega}{\tmu^{(0)}
(e_1, \ldots, \check{e}_l, \ldots, e_n)}
= 0,
\qquad (k=0).
\label{coef3}
\end{align}
First we absorb the $n$-form $\tmu^{(n)}$ to $H$ by defining 
$\tH = H + \rd \tmu^{(n)}$ in Equation \eqref{coef1}, and next,
using Equation \eqref{conditionofgsm21} step by step, Equations
\eqref{coef1}--\eqref{coef3} are rearranged to 
\begin{eqnarray}
&& \nabla \tmu^{(n-1)} = - (-1)^{n-1} \iota_{\rho} \tH,
\label{GNLSM1}
\\
&& \nabla \tmu^{(n-2)} + (-1)^{n-1}\, {}^E \rd^{\nabla} \tmu^{(n-1)}
= - (-1)^{(n-2)(n-1)}
\iota_{\rho}^{ 2} \tH,
\label{GNLSM2}
\\
&& \vdots
\nonumber
\\
&& \nabla \tmu^{(k-1)} + (-1)^{k}\, {}^E \rd^{\nabla} \tmu^{(k)}= - (-1)^{(k-1)\ldots (n-1)} \iota_{\rho}^{n+1-k} \tH,
\label{GNLSM3}
\\
&& \vdots
\nonumber
\\
&& {}^E \rd \tmu^{(0)} = - \iota_{\rho}^{n+1} \tH.
\label{GNLSM4}
\end{eqnarray}
Here the $E$-exterior covariant derivative ${}^E \rd^{\nabla}$ in Equations 
\eqref{GNLSM1}-\eqref{GNLSM4} is given by Eq.~\eqref{covariantLiederi}.
If we define
\begin{eqnarray}
&& 
\mu_k := (-1)^{(k+1) \ldots(n-1)} \tmu^{(k)},
\label{redefine1}
\end{eqnarray}
and
\begin{eqnarray}
&& \momega(\rho(e_1), \ldots, \rho(e_k), v_{k+1}, \ldots, v_{n+1})
:= 
(-1)^{n-k} \tH(\rho(e_k), \ldots, \rho(e_1), v_{k+1}, \ldots, v_{n+1}),
\label{redefine2}
\end{eqnarray}
Equations 
\eqref{GNLSM1}--\eqref{GNLSM4}
become Equations
\eqref{homotopyMS1}--\eqref{homotopyMS4}
in the definition of the homotopy momentum section.
Here $e_1, \ldots, e_k \in \Gamma(E)$ and 
$v_{k+1}, \ldots, v_{n+1} \in \mathfrak{X}(M)$.
Equation $\eqref{conditionofgsm21}$ is 
\begin{eqnarray}
\mu_{k-1} = \iota_{\rho} \mu_k,
\label{algecond}
\end{eqnarray}
for $k= 1, \ldots, n-1$.
We summarize our result. 
\begin{theorem}
We consider the $n$-dimensional nonlinear sigma model with $(n+1)$-dimensional WZ term. 
If we consider gauging by a Lie algebroid, 
the theory has a homotopy Hamiltonian Lie algebroid structure.
\end{theorem}
i.e. Suppose the gauge action functional $S_{LA}$ is Equation \eqref{ndgaugedsigmamodel}, the action functional $S_{LA}$ is gauge invariant if and only if $\mu_k := (-1)^{(k+1) \ldots(n-1)} \tmu^{(k)}$ are homotopy momentum sections with Equation \eqref{algecond}.

\section{Conclusion and discussion}
We have defined a homotopy momentum section on 
a Lie algebroid over a pre-multisymplectic manifold.
It is a simultaneous generalization of a momentum map on a symplectic manifold,
a homotopy momentum map on a pre-multisymplectic manifold,
and a momentum section on a symplectic manifold.

We analyzed relations between (weak) homotopy momentum sections and 
(weak) homotopy momentum maps on an action Lie algebroid. 
For it, the equivariant condition has been introduced.
If $\mu$ is equivariant, a homotopy momentum section is 
a homotopy momentum map.
On the other hand, a weak homotopy momentum section is always
a weak homotopy momentum map without the equivariance condition.

As an important application, we showed that the $n$-dimensional gauged nonlinear sigma model with WZ term with Lie algebroid gauging has a homotopy momentum section and a homotopy Hamiltonian Lie algebroid structure.

One of important applications of the momentum map theory is the equivalent cohomology and the symplectic reduction \cite{Marsden-Weinstein}.
These important applications are left for future works.
There is a discussion about reductions on a multisymplectic manifold
in the Lie group setting \cite{Blacker}.

The connection of a homotopy momentum section with the 
gauged nonlinear sigma model will provide possibility to analyze 
sigma models, field theories and their quantizations using the homotopy momentum section theory. Physical applications such as string theory can be discussed in this context.

\subsection*{Acknowledgments}
\noindent
The authors would like to thank 
Aliaksandr Hancharuk, Sylvain Lavau, Thomas Strobl, Alan Weinstein and Marco Zambon for useful discussion and comments.
This work was supported by the research promotion program for acquiring grants in-aid for scientific research(KAKENHI) in Ritsumeikan university
and JSPS KAKENHI Grant Number 22K03323.

\appendix
\section{Geometry of Lie algebroid}\label{geometryofLA}
We summarize notation, formulas and their local coordinate expressions of geometry of a Lie algebroid.

Let $(E, \rho, [-,-])$ be a Lie algebroid over a smooth manifold $M$.
$x^i$ is a local coordinate on $M$, $e_a \in \Gamma(E)$ is a basis of sections of $E$ and $e^a \in \Gamma(E^*)$ is a dual basis of sections of $E^*$. 
$i,j$, etc. are indices on $M$ and $a,b$, etc. are indices on the fiber of $E$.
Local coordinate expressions of the anchor map and the Lie bracket are
$\rho(e_a) f = \rho^i_a(x) \partial_i f$ and
$[e_a, e_b ] = C_{ab}^c(x) e_c$, where $f \in C^{\infty}(M)$ and $\partial_i = \tfrac{\partial}{\partial x^i}$.
Then, identities of $\rho$ and $C$ induced from the Lie algebroid condition are
\beqa 
&& \rho_a^j \partial_j \rho_{b}^i - \rho_b^j \partial_j \rho_{a}^i = C_{ab}^c \rho_c^i,
\label{LAidentity1}
\\
&& C_{ad}^e C_{bc}^d + \rho_a^i \partial_i C_{bc}^e + \mbox{Cycl}(abc) = 0.
\label{LAidentity2}
\eeqa

Let $\nabla$ be an ordinary connection on the vector bundle $E$.
The canonical \textit{$E$-connection} ${}^E \nabla: \Gamma(TM) \rightarrow \Gamma(TM \otimes E^*) $ on the space of vector fields $\mathfrak{X}(M) = \Gamma(TM)$ is defined by
\begin{eqnarray}
{}^E \nabla_{e} v &:=& \calL_{\rho(e)} v + \rho(\nabla_v e)
= [\rho(e), v] + \rho(\nabla_v e),
\end{eqnarray}
where 
$e \in \Gamma(E)$ and $v \in \Gamma(TM)$.
The standard $E$-connection on $E$ is 
\begin{eqnarray}
{}^E \nabla_{e} \alpha &:=& \nabla_{\rho(e)} \alpha.
\end{eqnarray}
for $\alpha \in \Gamma(\wedge^m E^*)$.
%
The dual $E$-connection ${}^E \nabla$ on the space of differential forms 
$\Omega^1(M)$ is obtained from the equation,
\begin{eqnarray}
{}^E \rd \bracket{v}{\alpha} = 
\bracket{{}^E \nabla e}{\alpha} + \bracket{e}{{}^E \nabla \alpha},
\end{eqnarray}
for a vector field $v$ and a $1$-form $\alpha$.
For a $1$-form $\alpha$, it is given by
\begin{eqnarray}
{}^E \nabla_{e} \alpha &:=& \calL_{\rho(e)} \alpha 
+ \bracket{\rho(\nabla e)}{\alpha}.
\end{eqnarray}
Let $\omega = \omega^b_{ai} \rd x^i \otimes e^a \otimes e_b$ be 
a connection $1$-form. 
For the basis vectors, the covariant derivatives are
${\nabla} e_a = - \omega_{ai}^b \rd x^i \otimes e_b$
and ${\nabla} e^a = \omega_{bi}^a \rd x^i \otimes e^b$.
Local coordinate expressions of 
covariant derivatives and the standard $E$-covariant derivative on $TM$ and $E$ are
\begin{eqnarray}
\nabla_i \alpha^a &=& \partial_i \alpha^a {+} \omega_{bi}^a \alpha^b,
\\
\nabla_i \beta_a &=& \partial_i \beta_a {-} \omega_{ai}^b \beta_b,
\\
{}^E \nabla_a v^i
&=&  \rho_a^j \partial_j v^i - \partial_j \rho^i_a v^j
{+} \rho^i_b \omega^b_{aj} v^j,
\\
{}^E \nabla_a \alpha_i
&=&  \rho_a^j \partial_j \alpha_i + \partial_i \rho^j_a \alpha_j
{-} \rho^j_b \omega^b_{ai} \alpha_j.
\end{eqnarray}

The covariant derivative of an $l$-form taking a value in 
$E^{\otimes m} \otimes E^{* \otimes n}$, 
$\alpha \in \Omega^l(M, E^{\otimes m} \otimes E^{* \otimes n})$ is
given by
\begin{align}
\nabla_j \alpha_{k_1 \ldots k_l}{}^{a_1 \ldots a_m}_{b_1 \ldots b_n} &=
\partial_j \alpha_{k_1 \ldots k_l}{}^{a_1 \ldots a_m}_{b_1 \ldots b_n} 
{+} \sum_{i=1}^m \omega_{cj}^{a_i} 
\alpha_{k_1 \ldots k_l}{}^{a_1 \ldots a_{i-1} c a_{i+1} \ldots a_m}_{b_1 \ldots b_n} 
{-} \sum_{i=1}^n \omega_{b_ij}^c 
\alpha_{k_1 \ldots k_l}{}^{a_1 \ldots a_m}_{b_1 \ldots b_{i-1} c b_{i+1} \ldots b_n}.
\end{align}

%
An $E$-torsion, a curvature, an $E$-curvature and a basic curvature,
$T \in \Gamma(E \oplus \wedge^2 E^*)$, $R \in \Omega^2(M, E \oplus E^*)$,
${}^E R \in \Gamma(\wedge^2 E^* \oplus E \oplus E^*)$,
and $S \in \Omega^1(M, \wedge^2 E^* \oplus E)$, are defined by
\beqa
R(s, s^{\prime}) &:=& [\nabla_s, \nabla_{s^{\prime}}] - \nabla_{[s, s^{\prime}]}, 
\\
T(s, s^{\prime}) &:=& {}^E\nabla_s s^{\prime} - {}^E\nabla_{s^{\prime}} s
- [s, s^{\prime}],
\label{Etorsion}
\\
{}^ER(s, s^{\prime}) &:=& [{}^E\nabla_s, {}^E\nabla_{s^{\prime}}] - {}^E\nabla_{[s, s^{\prime}]}, 
\\
S(s, s^{\prime}) &:=& \calL_s (\nabla s^{\prime}) 
- \calL_{s^{\prime}} (\nabla s) 
- \nabla_{\rho(\nabla s)} s^{\prime} + \nabla_{\rho(\nabla s^{\prime})} s
\nonumber \\ &&
- \nabla[s, s^{\prime}] = (\nabla T + 2 \mathrm{Alt} \, \iota_\rho R)(s, s^{\prime}),
\eeqa
for $s, s^{\prime} \in \Gamma(E)$.

\section{Supergeometry of Lie algebroid}\label{supergeometry}

Lie algebroids are described by means of $\mathbb{Z}$-graded geometry \cite{Vaintrob}. 
It is the so called $Q$-manifold description. 
A graded manifold $E[1]$ for a vector bundle $E$ are
shifted vector bundle spanned by local coordinates $x^i, \ (i=1, \ldots, \mathrm{dim}\, M)$ on the base manifold $M$ of degree zero, and local coordinates $q^a, \ (a=1, \ldots, \mathrm{rank}\, E)$ on the fiber of degree one, respectively. 
Degree one coordinate $q^a$ has the property, $q^a q^b = - q^b q^a$. $E$-differential forms which are 
sections of $\wedge^{\bullet} E^*$ are identified functions on the graded manifold $E[1]$, i.e., $C^{\infty}(E[1]) \simeq \Gamma(\wedge^{\bullet} E^*)$,
where the degree one odd coordinate $q^a$ is identified to the basis $e^a$ of sections of $E^*$.
The differential operator of degree $-1$, $\frac{\partial}{\partial q^a}$, is the derivation satisfying $\frac{\partial}{\partial q^a} q^b = \delta^b_a$, which is a linear operator on $C^{\infty}(E[1])$ satisfying the Leibniz rule.  

There exists a degree plus one (odd) vector field $Q$ on $E[1]$:
 \beq \label{Q}
 Q = \rho_a^i (x) q^a \frac{\partial}{\partial x^i} - \frac{1}{2} C^c_{ab}(x) q^a q^b  \frac{\partial}{\partial q^c} \, .
 \eeq
where $\rho_a^i (x)$ and $C^c_{ab}(x)$ are local functions of $x$. 
Then, the odd vector field $Q$ satisfies 
 \beq  Q^2 = 0 \, ,\label{Q2}
 \eeq
if and only if $\rho, C$ are the anchor map and the structure function of a Lie algebroid on $E$.
Identifying $C^{\infty}(E[1]) \simeq \Gamma(\wedge^\bullet E^*)$, 
$Q$ is the Lie algebroid differential ${}^E \rd$.

The precise correspondence of $Q$ with ${}^E\rd$ is as follows.
For $e^a$, the basis of $E^*$,
the map 
$j: \Gamma(\wedge^{\bullet}E^*) \rightarrow C^{\infty}(E[1])$
is given by the map of basis,
$j:(x^i, e^a) \mapsto (x^i, q^a)$.
The differential ${}^E\rd$ on $\Gamma(\wedge^{\bullet}E^*)$ is defined by the pullback using the map $j$, ${}^E\rd=j^*Q$.

Note that $Q$ is itself covariant even when no connection is introduced. In fact, introducing a connection we can rewrite $Q$ to the 
manifestly covariant form as
\begin{eqnarray}
Q &=& \rho^i_a(x) {q^a}
\left(\frac{\partial}{\partial x^i} - \omega_{bi}^c q^b \frac{\partial}{\partial q^c} \right)
+ \frac{1}{2} T^a_{bc}(x) q^b q^c \frac{\partial}{\partial q^a}
\nonumber \\
&=& \rho^i_a(x) q^a \nabla_i
+ \frac{1}{2} T^a_{bc}(x) q^b q^c \frac{\partial}{\partial q^a}.
\end{eqnarray}
Let $\alpha = \frac{1}{m!} \alpha_{a_1 \ldots a_m}(x) e^{a_1} \wedge \ldots \wedge e^{a_n} \in \Gamma(\wedge^m E^*)$ be an $E$-differential form. Then, the 
supergeometry description is $\ualpha := j_* \alpha = \frac{1}{m!} \alpha_{a_1 \ldots a_m}(x) q^{a_1} \ldots q^{a_n}$. The Lie algebroid differential 
${}^E \rd \alpha$ can be calculated by ${}^E \rd \alpha = j^* Q \ualpha$.

Next, we consider the Lie algebroid homology operator ${}^E \partial$.
Let $u = u^a(x) e_a \in \Gamma(E)$ be a section of $E$, where $e_a$
is a basis of the fiber. Then, the corresponding graded vector field $\uu$ is defined by $\uu = u^a (x) \partial/\partial q^a$.
Analogous to the Cartan calculus, 
the following \textit{derived bracket}
$[[\uu, Q], \uv]\, \ualpha$ for $u, v \in \Gamma(E)$ gives the ordinary Lie bracket taking the pullback to the nonsuper vector bundle $E$,
\begin{align}
\iota_{[u, v]} \alpha = j^* ([[\uu, Q], \uv]\, \ualpha).
\end{align}
Therefore, in the super formulation, the Lie algebroid homology operator is calculated as
\begin{align}
{}^E \partial_m (e_1 \wedge \ldots \wedge e_m)
&= j^* \sum_{1 \leq i < j \leq m} (-1)^{i+j} [[\ue_i, Q], \ue_j]\, \ue_1 \ldots \check{\ue_i} \ldots \check{\ue_j} \ldots \ue_m,
\label{defhomology1}
\end{align}
where 
$\ue_1 \ldots \check{\ue_i} \ldots \check{\ue_j} \ldots \ue_m$ in
the right hand side is super products of $\ue_i$.
For $f \in C^{\infty}(M)$, the formula of the Lie algebroid 
homology operator is
\begin{align}
{}^E \partial_m (f e_1 \wedge \ldots \wedge e_m)
&= f {}^E \partial_m (e_1 \wedge \ldots \wedge e_m)
+ \sum_{i=1}^m (-1)^{i-1} (\rho(e_i) f) (e_1 \wedge \ldots \wedge \check{e_i} \wedge \ldots \wedge e_m).
\end{align}
%

\newcommand{\bibit}{\sl}



\end{document}